\documentclass[11pt]{article}

\usepackage{amssymb,amsmath,latexsym}
\usepackage[dvips]{graphics}
\usepackage{mathrsfs}
\usepackage{extarrows}
\usepackage{float}
\usepackage{color}

\newtheorem{theorem}{Theorem}[section]

\numberwithin{equation}{section}

\def\ZZ{\mathbb{Z}}

\def\pf{\noindent{\bf Proof.} }
\def\qed{{\hfill $\Box$\medskip}}

\def\DDD{\mathcal{D}}

\def\SSS{\mathcal{S}}
\def\AAA{\mathcal{A}}
\def\BBB{\mathcal{B}}
\def\CCC{\mathcal{C}}
\def\TTT{\mathcal{T}}
\def\one{\uppercase\expandafter{\romannumeral1}}
\def\two{\uppercase\expandafter{\romannumeral2}}
\allowdisplaybreaks[4]

\input{epsf.sty}
\usepackage{epsfig}

\begin{document}
\title{\bf The Construction of Two Kinds of Bijections in Simple Random Walk Paths}
\author{Sai SONG\footnote{Key Laboratory of Advanced Theory and Application in Statistics and Data Science-MOE, School of Statistics, East China Normal University.}~~and Qiang YAO\footnote{Corresponding author. Key Laboratory of Advanced Theory and Application in Statistics and Data Science-MOE, School of Statistics, East China Normal University and NYU--ECNU Institute of Mathematical Sciences at NYU Shanghai.
                      E-mail: qyao@sfs.ecnu.edu.cn.}
                      }
\date{}
\maketitle{}

\begin{abstract}
It is known that for the $2n$--step symmetric simple random walk on $\ZZ$, two events have the same probability if and only if their sets of paths have the same cardinality. In this article, we construct two kinds of bijections between sets of paths with the same cardinality. The construction is natural and simple. It can be easily realized through programming. More importantly, this construction opens a door to prove that two events in the $2n$--step symmetric simple random walk on $\ZZ$ have the same probability and some further related results.
\end{abstract}

\noindent{\bf 2010 MR subject classification:}~05A19, 60C05


\section{Introduction}
We consider the $2n$--step simple random walk on the integer line $\ZZ$ starting from the origin, where $n$ is a fixed positive integer, which can be intuitively expressed as follows. At time $0$, there is a walker standing at the origin. Then at each step, he either moves to the left or to the right until he moves up to $2n$ steps. Rigorously,
denote $$\SSS_n:=\{(S_0,S_1,\dots,S_{2n}):~S_0=0,~|S_i-S_{i-1}|=1~\text{for any}~1\leq i\leq2n\}.$$ Then $\SSS_n$ is the collection of all $2n$--step simple random walk~(denoted by SRW for short) paths on $\ZZ$ starting from the origin. Therefore, $|\SSS_n|=2^{2n}$, where $|\cdot|$ denotes the cardinality of a set.\\

It is known that (1)~the set of $2n$--step SRW paths having endpoint at the origin~(denote this set by $\AAA_n$) has the same cardinality with the set of $2n$---step SRW paths not touching the origin except the starting point~(denote this set by $\BBB_n$); (2)~the set of $2n$--step SRW paths having endpoint at the origin and staying positive except the starting and ending points~(denote this set by $\CCC_n$) has the same cardinality with the set of $2n$---step SRW paths having endpoint at the origin and staying positive except the starting and ending points as well as only one other mid point~(denote this set by $\DDD_n$). The rigorous definitions of $\AAA_n$, $\BBB_n$, $\CCC_n$ and $\DDD_n$ will be given in Sections 2 and 3. Since two finite sets have the same cardinality if and only if bijections can be constructed between them, there are bijections between $\AAA_n$ and $\BBB_n$, as well as between $\CCC_n$ and $\DDD_n$. As far as we know, no literature provided the detailed construction before. In this article, we provide natural and simple bijections for the two facts respectively. The construction procedures of the two bijections can be easily realized computationally. The following Sections 2 and 3 are devoted to the facts (1) and (2) described above respectively. In Section 4 we will give some concluding remarks.

\section{Bijection \one}
Denote
\begin{align*}
\AAA_n&:=\{(0,S_1,\dots,S_{2n})\in\SSS_n:~S_{2n}=0\},\\
\BBB_n&:=\{(0,S_1,\dots,S_{2n})\in\SSS_n:~S_i\neq0~\text{for any}~1\leq i\leq2n\}.
\end{align*}
It is known that $|\AAA_n|=|\BBB_n|=C_{2n}^n$, where $C_n^m:=\dfrac{n!}{m!(n-m)!}~~(0\leq m\leq n)$ denotes the combinatorial number. The proof can be found in many textbooks, see for example, Lemma 4.3.3 of Durrett~(2010). By the symmetric property, to construct a bijection between $\AAA_n$ and $\BBB_n$, we only need to construct a bijection between $\AAA_n^\prime$ and $\BBB_n^\prime$, where
\begin{align*}
\AAA_n^\prime&:=\{(0,S_1,\dots,S_{2n})\in\SSS_n:~S_1>0,~S_{2n}=0\},\\
\BBB_n^\prime&:=\{(0,S_1,\dots,S_{2n})\in\SSS_n:~S_i>0~\text{for any}~1\leq i\leq2n\}.\\
\end{align*}
%The following result is well known. For completeness we give the proof here.
%\begin{prop}\label{p:result1}
%We have $|\AAA_n|=|\BBB_n|=\dfrac{1}{2}C_{2n}^n$.
%\end{prop}
%
%\pf $|\AAA_n|=\dfrac{1}{2}C_{2n}^n$ is obvious by symmetry. Let $\tau:=\inf\{k>0:~S_{2k}=0\}$. Then we have
%\begin{align}\label{e:B_nequiv}
%\BBB_n&=\{(S_1,\dots,S_{2n})\in\SSS_n:~S_1>0,~\tau>n\}\nonumber\\
%&=\{(S_1,\dots,S_{2n})\in\SSS_n:~S_1>0\}-\bigcup\limits_{k=1}^n\{(S_1,\dots,S_{2n})\in\SSS_n:~S_1>0,~\tau=k\},
%\end{align}
%where ``$-$'' denotes the difference of two sets when the latter one is a subset of the former one. By the reflection principle, we have for any $1\leq k\leq n$,
%\begin{align}\label{e:B_ne1}
%\left|\{(S_1,\dots,S_{2n})\in\SSS_n:~S_1>0,~\tau=k\}\right|&=(C_{2k-2}^{k-1}-C_{2k-2}^k)\cdot2^{2n-2k}\nonumber\\
%&=\frac{1}{2^{2k}}\cdot\frac{(2k-2)!}{k!(k-1)!}\cdot2^{2n}.
%\end{align}
%It is not difficult to prove by induction that
%\begin{equation}\label{e:B_ne2}
%\sum\limits_{k=1}^n\frac{1}{2^{2k}}\cdot\frac{(2k-2)!}{k!(k-1)!}=\frac{1}{2}-\frac{1}{2^{2n+1}}C_{2n}^n.
%\end{equation}
%Therefore, by (\ref{e:B_nequiv}), (\ref{e:B_ne1}) and (\ref{e:B_ne2}), we get
%\begin{align*}
%|\BBB_n|&=2^{2n-1}-2^{2n}\cdot\sum\limits_{k=1}^n\frac{1}{2^{2k}}\cdot\frac{(2k-2)!}{k!(k-1)!}\\
%&=2^{2n-1}-2^{2n}\cdot\left(\frac{1}{2}-\frac{1}{2^{2n+1}}C_{2n}^n\right)=\frac{1}{2}C_{2n}^n,
%\end{align*}
%as desired.\qed

In the following, we will construct a bijection between $\AAA_n^\prime$ and $\BBB_n^\prime$. For any $(0,S_1,\dots,S_{2n})\in\AAA_n^\prime$, denote $M:=\max\{S_1,\dots,S_{2n}\}$. Then define $a_M:=\min\{1\leq i\leq2n:~S_i=M\}$, $b_0:=0$, and inductively define $$a_k:=\min\{1\leq i\leq a_{k+1}:~S_i=k\},~~~~b_k:=\max\{1\leq i\leq a_{k+1}:~S_i=k\}$$ for $k=M-1,M-2,\dots,2,1$. Obviously, $$1=a_1\leq b_1<a_2\leq b_2<\dots<a_M$$ and $a_i=b_{i-1}+1$ for $i=1,\dots,M$.

\bigskip

For $l=1,\dots,2n$, define
\begin{equation}\label{e:phi}
T_\ell:=\begin{cases}2m-S_\ell~~~~\text{if}~a_m\leq\ell\leq b_m,~m=1,\dots,M-1,\\2M-S_\ell~~~~\text{if}~\ell\geq a_M.\end{cases}
\end{equation}
Next, define mapping $\Phi_1:~(0,S_1,\dots,S_{2n})\longmapsto(0,T_1,\dots,T_{2n})$. Since for any $a_m\leq\ell\leq b_m~(m=1,\dots,M-1)$, $$T_\ell=2m-S_\ell\geq m>0,$$ while for any $\ell\geq a_M$, $$T_\ell=2M-S_\ell\geq M>0,$$ we have shown that $\Phi_1$ is a mapping from $\AAA_n^\prime$ to $\BBB_n^\prime$.

\bigskip

\noindent\textbf{Remark.}~An intuitive expression of the mapping $\Phi_1$ is as follows. For any path in $\AAA_n$, we denote its maximum by $M~(>0)$ and divide its time period $[1,2n]$ by disjoint union $$[a_1,b_1]\cup\dots\cup[a_{M-1},b_{M-1}]\cup[a_M,2n].$$ For any $m=1,\cdots,M-1$, the path on the time period $[a_m,b_m]$ exhibits a ``valley like'' shape in the way that $S_{a_m}=S_{b_m}=m$ and $S_\ell<m$ for $\ell\in(a_m,b_m)$. Then the action of $\Phi_1$ is to reflect each ``valley'' at the height of the endpoint of its time interval and change it to a ``mountain like'' shape. The path in the rightmost time period is simply reflected at height $M$. Next, we give an explicit example for the case $n=3$. When $n=3$, we have $|\AAA_3^\prime|=|\BBB_3^\prime|=\dfrac{1}{2}C_6^3=10$. By the above construction, a bijection between $\AAA_3^\prime$ and $\BBB_3^\prime$ is
\begin{align*}
\Phi_1:~~~~~~~~~\AAA_3^\prime~~~~~~~~~&\longrightarrow~~~~~~~~~\BBB_3^\prime\\
(0,1,2,3,2,1,0)&\longmapsto(0,1,2,3,4,5,6)\\
(0,1,2,1,2,1,0)&\longmapsto(0,1,2,3,2,3,4)\\
(0,1,2,1,0,1,0)&\longmapsto(0,1,2,3,4,3,4)\\
(0,1,2,1,0,-1,0)&\longmapsto(0,1,2,3,4,5,4)\\
(0,1,0,1,2,1,0)&\longmapsto(0,1,2,1,2,3,4)\\
(0,1,0,1,0,1,0)&\longmapsto(0,1,2,1,2,1,2)\\
(0,1,0,1,0,-1,0)&\longmapsto(0,1,2,1,2,3,2)\\
(0,1,0,-1,0,1,0)&\longmapsto(0,1,2,3,2,1,2)\\
(0,1,0,-1,0,-1,0)&\longmapsto(0,1,2,3,2,3,2)\\
(0,1,0,-1,-2,-1,0)&\longmapsto(0,1,2,3,4,3,2)
\end{align*}
See Figure 1 for illustration. The blue line represents the original path in $\AAA_3^\prime$, and the red line represents its image in $\BBB_3^\prime$.

\begin{figure}[H]\label{picc:1}
 \center
 \includegraphics[width=10true cm, height=14true cm]{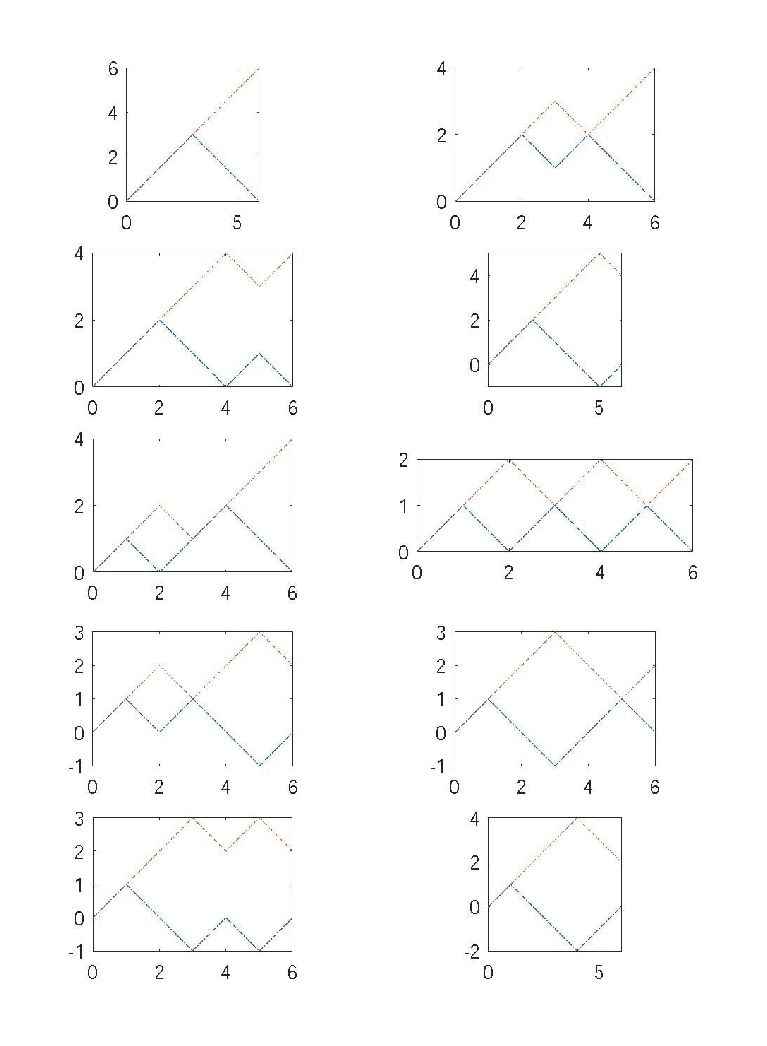}
 \caption{Illustration of $\Phi_1(\cdot)$ when $n=3$}
\end{figure}

\bigskip

\begin{theorem}\label{t:bijection1}
$\Phi_1$ is a bijection between $\AAA_n^\prime$ and $\BBB_n^\prime$.
\end{theorem}

\pf We can define the inverse mapping of $\Phi_1$ directly. For any $(0,T_1,\dots,T_{2n})\in\BBB_n^\prime$, denote $h:=\dfrac{1}{2}T_{2n}$. Then define $d_0:=0$, and inductively define $$c_k:=\min\{d_{k-1}\leq i\leq2n:~T_i=k\},~~~~d_k:=\max\{d_{k-1}\leq i\leq2n:~T_i=k\}$$ for $k=1,2,\dots,h$. Obviously, $$1=c_1\leq d_1<c_2\leq d_2<\dots<c_h\leq d_h$$ and $c_i=d_{i-1}+1$ for $i=1,\dots,h$.

\bigskip

For $\ell=1,\dots,2n$, define
\begin{equation}\label{e:psi}
S_\ell:=\begin{cases}2m-T_\ell~~~~\text{if}~c_m\leq\ell\leq d_m,~m=1,\dots,h-1,\\2h-T_\ell~~~~~\text{if}~\ell\geq c_h.\end{cases}
\end{equation}
Next, define the mapping $\Psi_1:~(0,T_1,\dots,T_{2n})\longmapsto(0,S_1,\dots,S_{2n})$. Since $S_{2n}=2h-T_{2n}=0$, it follows immediately that $\Psi_1$ is a mapping from $\BBB_n^\prime$ to $\AAA_n^\prime$.

\bigskip

We will show that $\Psi_1$ is the inverse mapping of $\Phi_1$ in the following. Comparing (\ref{e:phi}) with (\ref{e:psi})~(by taking $(0,T_1,\dots,T_{2n})=\Phi_1\left((0,S_1,\cdots,S_{2n})\right)$ where $(0,S_1,\dots,S_{2n})\in\AAA_n^\prime$), it can be seen that we only need to prove $h=M$ and $d_i=b_i$ for any $i=0,1,\dots,M-1$. Since $S_{2n}=0$, we have $$h=\frac{1}{2}T_{2n}=\frac{1}{2}(2M-S_{2n})=M.$$ Next, we will prove by induction that $d_i=b_i$ for any $i=0,1,\dots,M-1$. The fact $d_0=b_0$ is trivial. Suppose we have proved $d_{i-1}=b_{i-1}$ for some $i\in\{1,\dots,M-1\}$. We next prove $d_i=b_i$. On one hand, since $T_{d_{i-1}}=i-1$ and $T_{b_i}=i$, we have
\begin{equation}\label{e:dgeqb}
d_i=\max\{d_{i-1}\leq\ell\leq2n:~T_\ell=i\}\geq b_i.
\end{equation}
On the other hand, for any $\ell>b_i$, suppose $\ell\in[a_k,b_k]$ for $k>i$, then by the definition of $a_k$ and $b_k$, we have $S_\ell\leq k$. Therefore, $$T_\ell=2k-S_\ell\geq k>i.$$ Consequently, $$d_i=\max\{d_{i-1}\leq\ell\leq2n:~T_\ell=i\}\leq b_i.$$ Combined with (\ref{e:dgeqb}), we obtain the desired result $d_i=b_i$.\qed

\section{Bijection \two}

Denote $$\TTT_n:=\{(0,S_1,\dots,S_{2n})\in\SSS_n:~S_{2n}=0~\text{and}~S_i\geq0~\text{for any}~1\leq i\leq2n-1\}.$$
For any $(0,S_1,\dots,S_{2n})\in\TTT_n$, define $N\left((0,S_1,\dots,S_{2n})\right):=|\{i\in\{1,\dots,n-1\}:~S_{2i}=0\}|$. Then denote
\begin{align*}
\CCC_n&:=\{(0,S_1,\dots,S_{2n})\in\TTT_n:~N\left((0,S_1,\dots,S_{2n})\right)=0\},\\
\DDD_n&:=\{(0,S_1,\dots,S_{2n})\in\TTT_n:~N\left((0,S_1,\dots,S_{2n})\right)=1\}.
\end{align*}
By the reflection principle~(see for example, Theorem 4.3.1 of Durrett~(2010)), we have $$|\CCC_n|=C_{2n-2}^{n-1}-C_{2n-2}^n=\frac{(2n-2)!}{(n-1)!n!},$$ which is the $(n-1)$th Catalan number. The time homogeneity together with the well--known property of the Catalan numbers~(see for example, Theorem 3.3 in Roman~(2015)) lead to $$|\DDD_n|=\sum\limits_{k=1}^{n-1}|\CCC_k|\cdot|\CCC_{n-k}|=|\CCC_n|$$ for any $n\geq2$.

\bigskip

\noindent\textbf{Remark.}~We give an intuitive expression for the fact $|\CCC_n|=|\DDD_n|~(n\geq2)$ in the following. One day a gambler went to a casino without money and took fair games from which he earned~(paid) $1$ buck if he wined~(lost) each game, with debts forbidden. Suppose he took an alarm with him which could remind him when he was out of money. Then for any $n\geq2$, the probabilities of his alarm ringing for the first time and for the second time after the $2n$th game are the same.\\

%\begin{prop}\label{p:result2}
%For $n\geq2$, we have $|\CCC_n|=|\DDD_n|$.
%\end{prop}
%
%\pf By the reflection principle, we have $$|\CCC_n|=C_{2n-2}^{n-1}-C_{2n-2}^n=\frac{(2n-2)!}{(n-1)!n!}.$$
%Then by the time homogeneity, we have for $n\geq2$, $$|\DDD_n|=\sum\limits_{k=1}^{n-1}|\CCC_k|\cdot|\CCC_{n-k}|.$$ Denote $a_n:=\dfrac{(2n-2)!}{(n-1)!n!}$ for $n\geq1$, then we need to prove
%\begin{equation}\label{e:combinduction}
%a_n=\sum\limits_{k=1}^{n-1}a_ka_{n-k}
%\end{equation}
%for any $n\geq2$. Define
%\begin{equation}\label{e:fexpansion1}
%f(t):=\sum\limits_{n=1}^\infty a_nt^n.
%\end{equation}
%Then by the power series expansion, we have $f(t)=\dfrac{1-\sqrt{1-4t}}{2}$ for $t\leq\dfrac{1}{4}$. Therefore, we have $$[f(t)]^2-f(t)+t=0.$$ Furthermore, by Fubini's Theorem, we have
%\begin{equation}\label{e:fexpansion2}
%f(t)=t+[f(t)]^2=t+\sum\limits_{k=1}^\infty a_kt^k\sum\limits_{n=k+1}^\infty a_{n-k}t^{n-k}=t+\sum\limits_{n=2}^\infty\left(\sum\limits_{k=1}^{n-1}a_ka_{n-k}\right)t^n.
%\end{equation}
%Comparing (\ref{e:fexpansion1}) with (\ref{e:fexpansion2}), we can get from the uniqueness of the power series expansion of $f$ that $a_n=\sum\limits_{k=1}^{n-1}a_ka_{n-k}$ for any $n\geq2$. Therefore, (\ref{e:combinduction}) is true, as desired.\qed

Now we construct a bijection between $\CCC_n$ and $\DDD_n$, where $n\geq2$. For any $(0,S_1,\dots,S_{2n})\in\CCC_n$, denote $\tau:=\min\{k>1:~S_k=1\}$. Then by the property of $\CCC_n$, we have $2<\tau<2n$.

\bigskip

For $\ell=1,\dots,2n$, define
\begin{equation}\label{e:phi2}
T_\ell:=\begin{cases}S_\ell-2~~~~\text{if}~1<\ell<\tau~\text{and}~S_{\ell+1}=S_\ell-1,\\~~~S_\ell~~~~~~\text{otherwise}.\end{cases}
\end{equation}
Next, define mapping $\Phi_2:~(0,S_1,\dots,S_{2n})\longmapsto(0,T_1,\dots,T_{2n})$. For any $(0,S_1,\dots,S_{2n})\in\CCC_n$, we have $S_{2n}=0$ and $S_\ell>0$ for any $1\leq\ell\leq2n-1$. Then by (\ref{e:phi2}), we have $$T_{2n}=S_{2n}=0,~~T_{\tau-1}=S_{\tau-1}-2=0,~~T_\ell=S_\ell>0~\text{for any}~\tau\leq\ell\leq2n-1.$$ Furthermore, for $1\leq\ell\leq\tau-2$, if $S_{\ell+1}=S_\ell+1$, then $T_\ell=S_\ell>0$, while if $S_{\ell+1}=S_\ell-1$, then $S_\ell=S_{\ell+1}+1>2$ and therefore $T_\ell=S_\ell-2>0$. So we have $(0,T_1,\dots,T_{2n})\in\DDD_n$, and therefore, $\Phi_2$ is a mapping from $\CCC_n$ to $\DDD_n$.

\bigskip

\noindent\textbf{Remark.}~An intuitive expression of the mapping $\Phi_2$ is as follows. For any path in $\CCC_n$, a point among the path drops by $2$ if its time spot is in the time period $[2,\tau-1]$ and it descends immediately afterwards, where $\tau$ is the first time spot that the path touches $1$ in the time period $[3,2n-1]$. Otherwise, the point remains in its intrinsic position. Next, we give an explicit example for the case $n=4$. When $n=4$, we have $|\CCC_4|=|\DDD_4|=\dfrac{6!}{3!\cdot4!}=5$. By the above construction, a bijection between $\CCC_4$ and $\DDD_4$ is
\begin{align*}
\Phi_2:~~~~~~~~~~~~~\CCC_4~~~~~~~~~~~~&\longrightarrow~~~~~~~~~~~~\DDD_4\\
(0,1,2,3,4,3,2,1,0)&\longmapsto(0,1,2,3,2,1,0,1,0)\\
(0,1,2,3,2,3,2,1,0)&\longmapsto(0,1,2,1,2,1,0,1,0)\\
(0,1,2,3,2,1,2,1,0)&\longmapsto(0,1,2,1,0,1,2,1,0)\\
(0,1,2,1,2,3,2,1,0)&\longmapsto(0,1,0,1,2,3,2,1,0)\\
(0,1,2,1,2,1,2,1,0)&\longmapsto(0,1,0,1,2,1,2,1,0)
\end{align*}
See Figure 2 for illustration. The blue line represents the original path in $\CCC_4$, and the red line represents its image in $\DDD_4$.

\begin{figure}[H]\label{picc:2}
 \center
 \includegraphics[width=10true cm, height=8true cm]{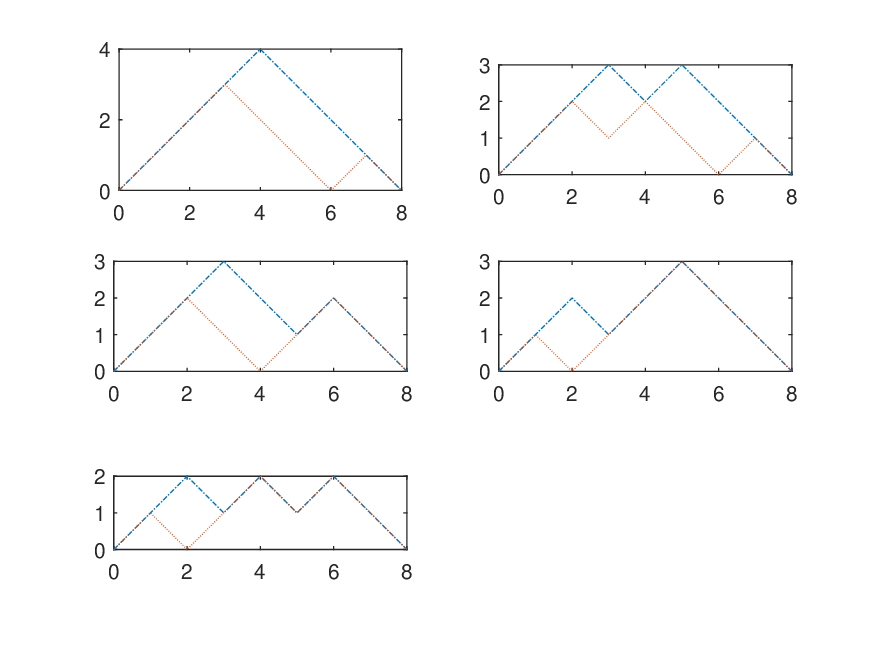}
 \caption{Illustration of $\Phi_2(\cdot)$ when $n=4$}
 \end{figure}

\bigskip

\begin{theorem}\label{t:bijection2}
$\Phi_2$ is a bijection between $\CCC_n$ and $\DDD_n$.
\end{theorem}

\pf We can define the inverse mapping of $\Phi_2$ directly. For any $(0,T_1,\dots,T_{2n})\in\DDD_n$, denote $\nu:=\min\{k>0:~T_k=0\}$. Then by the property of $\DDD_n$, we have $1<\nu<2n-1$.

\bigskip

For $\ell=1,\dots,2n$, define
\begin{equation}\label{e:psi2}
S_\ell:=\begin{cases}T_\ell+2~~~~\text{if}~1<\ell\leq\nu~\text{and}~T_\ell=T_{\ell-1}-1,\\~~~T_\ell~~~~~~\text{otherwise}.\end{cases}
\end{equation}
Next, define mapping $\Psi_2:~(0,T_1,\dots,T_{2n})\longmapsto(0,S_1,\dots,S_{2n})$. For any $(0,T_1,\dots,T_{2n})\in\DDD_n$, we have $T_\nu=T_{2n}=0$, and $T_\ell>0$ for any $\ell\in\{1,\dots,2n-1\}\setminus\{\nu\}$. Then by (\ref{e:psi2}), when $\ell\in\{1,\dots,2n-1\}\setminus\{\nu\}$, we have $$S_\ell\geq T_\ell>0,$$ while $S_\nu=T_\nu+2>0$ and $S_{2n}=T_{2n}=0$. So we have $(0,S_1,\dots,S_{2n})\in\CCC_n$, and therefore, $\Psi_2$ is a mapping from $\DDD_n$ to $\CCC_n$.

\bigskip

We will show that $\Psi_2$ is the inverse mapping of $\Phi_2$ in the following. Comparing (\ref{e:phi2}) with (\ref{e:psi2})~(by taking $(0,T_1,\dots,T_{2n})=\Phi_2\left((0,S_1,\dots,S_{2n})\right)$ where $(0,S_1,\dots,S_{2n})\in\CCC_n$), it can be seen that we only need to prove $\nu=\tau-1$, and $S_{\ell+1}=S_\ell-1$ is equivalent to $T_\ell=T_{\ell-1}-1$ for any $1<\ell<\tau$. Since $1<\tau-1<\tau$ and $S_\tau=S_{\tau-1}-1$, we can get from (\ref{e:phi2}) that $T_{\tau-1}=0$. Therefore, by the uniqueness of $\ell\in\{1,\cdots,2n-1\}$ such that $T_\ell=0$, we get $\nu=\tau-1$. Next, for any $1<\ell<\tau-1$, if $S_{\ell+1}=S_\ell-1$, then $$T_\ell-T_{\ell-1}=\begin{cases}(S_\ell-2)-S_{\ell-1}=-1~~~~~~~~~~~~\text{if}~S_\ell=S_{\ell-1}+1,\\(S_\ell-2)-(S_{\ell-1}-2)=-1~~~~~\text{if}~S_\ell=S_{\ell-1}-1.\end{cases}$$
On the other hand, if $T_\ell=T_{\ell-1}-1$, then $$S_{\ell+1}-S_\ell=\begin{cases}T_{\ell+1}-(T_\ell+2)=-1~~~~~~~~~~~~~~\text{if}~T_{\ell+1}=T_\ell+1,\\(T_{\ell+1}+2)-(T_\ell+2)=-1~~~~~~\text{if}~T_{\ell+1}=T_\ell-1.\end{cases}$$
The proof of Theorem \ref{t:bijection2} is complete.\qed

\section{Concluding remarks}

In this article, we provide two explicit bijections between sets of paths with the same cardinality in $2n$--step simple random walks on $\ZZ$. The two projections can be easily realized computationally. For example, Figure 1 and Figure 2 above are both drawn by the construction procedure of the two bijections using MATLAB. In addition, the construction of these bijections opens a door to prove some related results. For example, the above Bijection \two~ can provide a new probabilistic way to prove the explicit form of the Catalan numbers given the recurrence formula and the initial value. The idea of the proof procedure is as follows. By the construction of Bijection \two, we know that the $(n-1)$th Catalan number is just the number of paths in the $2n$--step simple random walk that start and end at the origin and stay positive in the middle. Together with the reflection principle, the explicit form of the $(n-1)$th Catalan number can be obtained. We remark here that this probabilistic proof procedure is beyond the $6$ proof methods exhibited on Wikipedia. We hope this article could inspire the construction of other bijections in the future.

\bigskip

%\noindent\textbf{Statement.}\quad  On behalf of all authors, the corresponding author states that there is no conflict of interest.

\noindent\textbf{Acknowledgment.}\quad The research of Qiang Yao was partially supported by the Natural Science Foundation of China~(No.11671145), the program of China Scholarships Council~(No.201806145024), and the 111 Project~(B14019).

\end{document}